\documentclass{amsart}

\newcommand{\PP}{{\mathbb P}}
\newcommand{\bd}{\mathbf{d}}
\newcommand{\diag}{\operatorname{diag}\nolimits}
\newcommand{\im}{\operatorname{im}}
\newcommand{\la}{\langle}
\newcommand{\liea}[1]{\mathfrak{#1}}
\newcommand{\lieg}[1]{\mathrm{#1}}
\newcommand{\ra}{\rangle}
\newcommand{\rk}{\operatorname{rk}}
\newcommand{\tr}{\operatorname{tr}\nolimits}

\theoremstyle{plain}
\newtheorem{thm}{Theorem}[section]
\newtheorem{lm}[thm]{Lemma}
\newtheorem{cor}[thm]{Corollary}
\newtheorem{prop}[thm]{Proposition}

\theoremstyle{definition}
\newtheorem{re}[thm]{Remark}

\begin{document}

\title{Higher secant varieties of the minimal adjoint orbit}
\author{Karin Baur and Jan Draisma}
\thanks{Both authors are supported by the Swiss National Science Foundation}
\address{Karin Bauer, Department of Mathematics, University of California, 
	San Diego, USA}
\email{kbaur@math.ucsd.edu}
\address{Jan Draisma, Mathematisches Institut Universit\"at Basel, Switzerland}
\email{jan.draisma@unibas.ch}
\date{16 December 2003}

\maketitle

\section*{Abstract}

The adjoint group of a simple complex Lie algebra $\liea{g}$ has a unique
minimal orbit in the projective space $\PP\liea{g}$, whose pre-image
in $\liea{g}$ we denote by $C$. We explicitly describe, for every
classical $\liea{g}$ and every natural number $k$, the Zariski closure
$\overline{kC}$ of the union $kC$ of all spaces spanned by $k$ points
on $C$. The image of this set in $\PP \liea{g}$ is usually called the
$(k-1)$-st {\em secant variety} of $\PP C$, and its dimension and defect
are easily determined from our explicit description. In particular,
it follows that the smallest $k$ for which $\overline{kC}$ is equal
to $\liea{g}$, is $n$ for $\liea{sl}_n$, $2n$ for $\liea{sp}_{2n}$, $4$
for $\liea{o}_7$, and $\lfloor \frac{n}{2} \rfloor$ for $\liea{o}_n,\
n \geq 8$; we find that the upper bound on this $k$ provided by a theorem
of Zak on secants of general varieties, is off by a factor of $2$ in the
cases of $\liea{sl}_n$ and $\liea{o}_n$, but sharp for $\liea{sp}_{2n}$.

The orthogonal Lie algebras turn out to be the most difficult, by far:
while all sets $kC$ are closed in the other two cases, this is not
true for $2C$ in $\liea{o}_n$, and we discuss the problems arising in
describing the sets $kC$ for $k \geq 3$. In particular, we do not know
the smallest $k$ for which $kC$ is equal to $\liea{o}_n$, though we do
prove that it is at most $\lfloor \frac{n}{2} \rfloor + 3$.

\section{Introduction and results}

The projective space $\PP \liea{g}$, where $\liea{g}$ is a semisimple
Lie algebra over an algebraically closed field $K$ of characteristic
zero, has a unique minimal orbit under the action of the adjoint group
of $\liea{g}$; let $C$ denote the pre-image of this orbit in $\liea{g}$.
Here `minimal' refers to the inclusion order among orbit closures, so the
minimality of $\PP C$ means that it is contained in the closure of any
other orbit. The set $C$, itself a nilpotent orbit, plays an important
role in several branches of Lie theory: First, $C$ consists of all
long root vectors relative to appropriate Cartan subalgebras (or of all
highest root vectors relative to Borel subalgebras) and is therefore of
interest in representation theory. Alternatively, $C$ may be described
as the set of all non-zero $X \in \liea{g}$ for which $[X,[X,\liea{g}]]
\subseteq KX$ \cite{Kaji99}, and these {\em extremal elements} pop
up in the classification of Lie algebras in positive characteristic
\cite{Cohen2001,Premet97} (for a possible connection between our
results and those of \cite{Cohen2001}, see the conclusion of this paper).
We, now, are to discuss properties of $C$ that are interesting from a
geometric point of view, namely: what do the higher secant varieties of
$\PP C$ in $\PP \liea{g}$ look like, and what are the corresponding
defects of $\PP C$? This work is part of a larger project, which
asks for the secant varieties of the minimal orbit in any irreducible
representation of any reductive algebraic group.

Recall that the $(k-1)$-st {\em secant variety} of $\PP C$ in the
projective space $\PP \liea{g}$ is the Zariski closure of the union
of all projective subspaces of $\PP \liea{g}$ spanned by $k$ points on
$\PP C$. As $C$ is a cone, this secant variety is in fact equal to
$\PP(\overline{kC} \setminus \{0\})=\overline{\PP(kC\setminus\{0\})}$
where $kC$ is the set defined by
\[ kC:=\{J_1+\ldots+J_k \mid J_i \in C \text{ for all } i\}. \]
The {\em expected dimension} of $kC$ is $\min \{k \dim C, \dim
\liea{g}\}$, and this expected dimension minus the actual dimension
of $kC$ is usually called the {\em $(k-1)$-defect} of $C$ (or of
$\PP C$). 

We are to present explicit descriptions of the sets $\overline{kC}$
in the case where $\liea{g}$ is classical, which extend the results
on the {\em first} secant variety of $\PP C$ by Kaji {\em et al}
\cite{Kaji99,Kaji2000}. It should be mentioned that their method applies
to the exceptional simple Lie algebras, as well, while it is not obvious
how to uniformise our case-by-case approach so as to incorporate those
in our treatment. Closely related to the matter of this paper is also
\cite{Catalisano2002}, which treats the higher secant varieties of the
variety of indecomposable tensors in a tensor product.

The research on higher secant varieties of general varieties
finds its origin in the works of Palatini and Terracini
\cite{Palatini09,Terracini11}, and an important part of this research
concerns (bounds on) the dimensions of secant varieties, as well as
the construction of concrete varieties attaining these dimensions
\cite{Adlandsvik87,Adlandsvik88,Catalano96,Fujita81,Goldstein85,
Landsberg96,Zak93}. The monograph \cite{Zak93} by Zak contains
the following result that we compared to our concrete situation: suppose
that the first defect $\delta$ of an irreducible smooth $n$-dimensional
projective algebraic variety $X$, lying in and spanning $\PP^N$, is
non-zero. Then the smallest $k$ for which the $k$-th secant variety of $X$
is equal to $\PP^N$, is at most $\lfloor\frac{n}{\delta}\rfloor$. Though
derived for application to varieties that---unlike the minimal orbit
$\PP C$---have low codimension in the ambient projective space, this
bound turns out to be quite good for the minimal orbit: it is roughly
twice the actual value for $\liea{sl}_n$ and $\liea{o}_n$, and it is
sharp for $\liea{sp}_{2n}$.

Turning our attention to a classical Lie algebra $\liea{g}$, we define
the {\em rank} of an element $A \in \liea{g}$, denoted $\rk(A)$, to
be its rank as a linear map in the standard $\liea{g}$-module $V$. In
the cases of $\liea{sl}_n$ and $\liea{sp}_{2n}$ the minimal orbit $C$
consists of all elements of rank $1$, and the following theorem, the
subject of Section \ref{sec:slsp}, identifies the secant varieties of
study as certain determinantal varieties.

\begin{thm} \label{thm:slsp}
If $\liea{g}=\liea{sl}_n\ (n \geq 2)$ or $\liea{g}=\liea{sp}_{2n}\
(n \geq 2)$, then we have $kC=\overline{kC}=\{A \in \liea{g} \mid \rk(A)
\leq k\}$ for all $k,\ 2 \leq k \leq n$.
\end{thm}

This is not new for $\liea{sl}_n$: \cite{Kraft2001} contains a proof
using the Jordan normal form, while our proof was inspired by \cite[\S 56,
Exercise 6]{Halmos74}.

\begin{cor} \label{cor:slsp}
\begin{enumerate}
\item For $\liea{g}=\liea{sl}_n$ and $1 \leq k \leq n$ the dimension of
$kC$ is $2kn-k^2-1$, so that the $(k-1)$-defect of $C$ is $\min \{(k-1)^2,
(n-k)^2\}$.

\item For $\liea{g}=\liea{sp}_{2n}$ and $1 \leq k \leq 2n$ the
dimension of $kC$ is $\binom{2n+1}{2} - \binom{2n+1-k}{2}$, so that the
$(k-1)$-defect of $C$ is $\min\{\binom{k}{2},\binom{2n+1-k}{2}\}$.
\end{enumerate}
\end{cor}

The result for $\liea{o}_n$ is radically different; we assume $n \geq 7$
here, as the other (simple) cases are dealt with by the preceding 
theorem. Now we have 
\[ C=\{ J \in \liea{o}_n \mid \rk(J)=2 \text{ and } J^2=0\} \text{
	(see Subsection \ref{ssec:C})}, \]
and one might hope that $kC$ is simply the set of all elements of rank at
most $2k$---but this is not true! To describe the first and the second
secant variety, denote by $S_2 \subseteq \liea{o}_n$ the set of all
semisimple elements of rank $4$ whose non-zero eigenvalues (on $V$)
are $a,a,-a,-a$ for some $a \in K^*$. Similarly, let $S_3 \subset
\liea{o}_n$ be the set of all semisimple elements of $\liea{o}_n$ of
rank $6$ with $6$ distinct non-zero eigenvalues $a,b,c,-a,-b,-c \in
K$ satisfying $a+b+c=0$. We then have the following theorem (Section
\ref{sec:o}).

\begin{thm} \label{thm:o}
For $\liea{g}=\liea{o}_n\ (n \geq 7)$ the sets $\overline{2C}$ and
$\overline{3C}$ are equal to $\overline{S_2}$ and $\overline{S_3}$,
respectively, while for $k \geq 4$ we have $\overline{kC}=\{A \in
\liea{o}_n \mid \rk(A) \leq 2k\}$.
\end{thm}

\begin{cor} \label{cor:o}
For $\liea{g}=\liea{o}_n\ (n \geq 7)$ the dimensions of $C, 2C, 3C$,
and $kC\ (4 \leq k \leq \lfloor \frac{n}{2} \rfloor)$ are $2n-6, 4n-13,
6n-22$, and $\binom{n}{2}-\binom{n-2k}{2}$, respectively (the dimension
of $4C$ in $\liea{o}_7$ is $\binom{7}{2}$).  Hence, the $(k-1)$-defect of
$C$ is equal to $1$ for $k=2$, equal to $4$ for $k=3$, and equal to $\min
\{k(2k-5), \binom{n-2k}{2}\}$ if $k \geq 4$ (and zero for $(k,n)=(4,7)$).
\end{cor}

Note that Theorem \ref{thm:o} only mentions the closures $\overline{kC}$,
not the sets $kC$ themselves. This is because we do not know their
exact structure; we now list what we do know.  First, the set $2C$
is already not closed; indeed, Kaji {\em et al} determined the
nilpotent orbits lying in its closure \cite{Kaji2000}, and it
turns out that $\overline{2C} \setminus 2C$ consists of a single
such orbit. To formulate our proposition to that effect, recall
that nilpotent orbits of $\lieg{O}_n$ on $\liea{o}_n$ correspond,
through the Jordan normal form, to partitions of $n$ whose even
entries have even multiplicities. If $\bd=(d_1,\ldots,d_m),\ d_1
\geq d_2 \geq \ldots \geq d_m$ is such a partition, then we denote by
$\mathcal{O}[\bd]=\mathcal{O}[d_1,\ldots,d_m]$ the corresponding nilpotent
orbit. For example, in this notation we have $C=\mathcal{O}[2,2,1^{n-4}]$.

\begin{prop} \label{prop:2C}
The set $2C$ is equal to $\overline{2C} \setminus \mathcal{O}[3,2,2,1^{n-7}]$.
\end{prop}

The fact that $\overline{2C} \setminus 2C$ is a nilpotent orbit suggests
to determine, for a general nilpotent orbit $\mathcal{O}$, the smallest
$k$ for which $\mathcal{O}$ is contained in $kC$. Our partial result in
this direction uses the notation $l(\bd):=|\{i\mid d_i\ \mbox{is odd},\
d_i>1\}|$. Furthermore, by the {\em rank} of an orbit $\mathcal{O}$
we shall mean the rank of an element of that orbit.

\begin{thm} \label{thm:nilp}
Let $\bd$ be a partition of $n$ as above, and let $2k$ be the rank of
$\mathcal{O}[\bd]$. Then $\mathcal{O}[\bd]$ is contained in $(k+1)C$. If
moreover $l(\bd)$ is even, or if $l(\bd)$ is odd and $d_1>5$, then
$\mathcal{O}(\bd)$ is already contained in $kC$.
\end{thm}

The upper bound $k+1$ (notation as in the preceding theorem) is sharp
for $\bd=[3,2,2,1^{n-7}], [3,1^{n-3}]$, and $[5,1^{n-5}]$.  Hence,
the nilpotent orbits of smallest rank for which we do not know the
smallest $kC$ containing them, are $\mathcal{O}[3,3,3,1^{n-9}]$ and
$\mathcal{O}[5,2,2,1^{n-9}]$, both of rank $6$. In conclusion, it seems
hard to write a general element of $\liea{o}_n$ as a sum of as few as
possible elements of $C$. The rank reduction argument used to prove
Theorem \ref{thm:o}, however, does give an upper bound to the maximum
number of terms needed.

\begin{thm} \label{thm:o2}
Every element of $\liea{o}_n$ having rank at most $2k$ lies in
$(k+3)C$. In particular, $(\lfloor \frac{n}{2} \rfloor+3)C=\liea{o}_n$.
\end{thm}

\subsection*{Acknowledgments}
We thank Hanspeter Kraft for motivating discussions on the subject of this
paper, and Jochen Kuttler for his short proof of Lemma \ref{lm:BilForms}
below.

\section{Appetisers: $\liea{sl}_n$ and $\liea{sp}_n$} \label{sec:slsp}

For all classical simple Lie algebras $\liea{g}$ the minimal orbit $C$
consists of matrices of some low rank $r$ ($r=1$ for $\liea{sl}_n$
and $\liea{sp}_{2n}$, and $r=2$ for $\liea{o}_n$; see below). As a
result, an element of $kC$ has rank at most $kr$. Now if $A$ lies in
$kC$, then by definition there exists a $J \in C$ such that $A-J \in
(k-1)C$, hence if $A$ has the maximal possible rank $kr$, then its rank
must decrease by $r$ upon subtracting $J$: $\rk(A-J)=\rk(A)-r=(k-1)r$.
It seems therefore reasonable, given an element $A \in \liea{g}$ that
we want to write as a sum of points on $C$, to look for $J \in C$ such
that $\rk(A)-\rk(J) < \rk(A)$. The easy lemma below turns out to be an
effective tool in the quest for such $J$.

Here, and in the rest of this paper, $V$ stands for the standard module
of the classical Lie algebra under consideration, $V^*$ denotes the
linear dual of $V$, and $\langle .,.  \rangle$ is the natural pairing
$V \times V^* \rightarrow K$. We identify $\liea{gl}(V)$ with $V \otimes
V^*$, and recall that under this identification the rank one elements of
$\liea{gl}(V)$ correspond to the tensors $y \otimes \eta$ with non-zero
$y \in V$ and $\eta \in V^*$. Furthermore, for $A \in \liea{gl}(V)$
we define the dual map $A^* \in \liea{gl}(V^*)$ by $\langle x,A^*\xi
\rangle=\langle Ax,\xi \rangle$.

\begin{lm} \label{lm:Rank1}
For $A \in \liea{gl}(V)$ and non-zero elements $y \in V,\eta \in V^*$ we have
\begin{enumerate}
\item $\rk(A-y \otimes \eta)<\rk(A)$ if and only if $y \in \im A, \ker
\eta \supseteq \ker A$, and $\langle x, \eta \rangle=1$ for some (and
hence for any) $x \in A^{-1}y$; and
\item $y \otimes \eta$ lies in $\liea{sl}(V)$ if and only if $\langle y, \eta
\rangle=0$.

\end{enumerate}
\end{lm}

\begin{re}
The conditions in the first statement are self-dual in $y$ and $\eta$,
and can be rephrased more compactly as $\langle A^{-1} y, \eta
\rangle=\{1\}$ or as $\langle y,(A^*)^{-1} \eta \rangle=\{1\}$.
\end{re}

\begin{proof}
The second statement is obvious. As for the first statement: if $y$ and
$\eta$ satisfy the stated conditions, then $\ker(A-J)=\ker A \oplus
Kx$, so that $\rk(A-J)<\rk(A)$ as claimed. Conversely, suppose that
$\rk(A-J)<\rk(A)$, and let $x'$ be an element of $\ker(A-J) \setminus
\ker A$. Then $(A-J)x'=Ax' - \langle x', \eta \rangle y=0$ while $Ax'
\neq 0$. We conclude that $c:=\langle x', \eta \rangle$ is non-zero,
so that we may set $x:=x'/c$. Now $\langle x,\eta \rangle=1$ and $Ax=y$
and, by a dual argument, $\eta$ lies in $\im A^*$, which is equivalent
to $\ker \eta \supseteq \ker A$.
\end{proof}

Now first consider the Lie algebra $\liea{sl}_n$ with standard module
$V=K^n$. The highest root vector with respect to the usual Cartan and
Borel subalgebras of $\liea{sl}_n$ is the matrix with zeroes everywhere
and a $1$ in the upper right corner, hence of rank one. It is easy to
see that the rank one elements of $\liea{sl}_n$ form one orbit under
the adjoint group $\lieg{PSL}_n$: the minimal orbit $C$. Now we are
ready to prove Theorem \ref{thm:slsp} in the case where
$\liea{g}=\liea{sl}_n$.

\begin{proof}[Proof of Theorem \ref{thm:slsp} for $\liea{g}=\liea{sl}_n$.]
Let $A \in \liea{sl}_n$ be non-zero; we show that there exists a $J \in C$
for which $\rk(A-J)=\rk(A)-1$. Indeed, $A$ induces a linear map $\bar{A}$
on the non-zero space $V/\ker A$, and we have $\tr \bar{A}=0$.  As $K$
has characteristic $0$, the map $\bar{A}$, having trace zero, cannot be a
non-zero scalar, so that there exists an $x \in V$ for which $x+\ker A$
is not a scalar multiple of $\bar{A}(x+\ker A)$.  This means that $x$
does not lie in $KAx+\ker A$, hence there exists a linear function $\eta
\in V^*$ that vanishes on $KAx+\ker A$ but has the value $1$ on $x$. Now
$J:= Ax \otimes \eta$ lies in $\liea{sl}_n$ and has $\rk(A-J)=\rk(A)-1$
by Lemma \ref{lm:Rank1}. 

By induction, this proves that every element of $\liea{sl}_n$ of rank $k$
lies in $kC$, for $k=1,\ldots,n$. The inclusions $kC \subseteq (k+1)C$
and $0 \in 2C$, following from the fact that $C$ is a cone, conclude the
proof that for $k \geq 2$ every element of rank {\em at most} $k$ lies
in $kC$. Conversely, it was observed at the beginning of this section
that $kC$ does not contain elements of rank higher than $k$.
\end{proof}

We proceed to prove Theorem \ref{thm:slsp} for the symplectic Lie algebra
$\liea{sp}_{2n}$. Let $V=K^{2n}$ be the standard $\liea{sp}_{2n}$-module,
and denote by $(.,.)$ the non-degenerate skew bilinear form on $V$ defining
$\liea{sp}_{2n}$. Define the linear maps $\phi:V \rightarrow V^*$
and, for $A \in \liea{gl}(V)$, $A^T:V \rightarrow V$ by $\langle x,
\phi y \rangle=(x,y)$ and $(x,A^T y)=(Ax,y)$ for $x,y \in V$ and $\eta
\in V^*$. We first describe the minimal $\lieg{Sp}_{2n}$-orbit $C$
in $\liea{sp}_{2n}$ in a convenient way.

\begin{lm} \label{lm:Csp2n}
$C=\{ x \otimes \phi x \mid x \in V\setminus\{0\} \} = \{ J \in
	\liea{sp}_{2n} \mid \rk J =1\}.$
\end{lm}

\begin{proof}
First, the highest root vector in $\liea{sp}_{2n}$ (with respect to some
choice of Cartan and Borel subalgebras) is easily seen to have rank one.
Now a rank one element $J=y \otimes \eta$ of $\liea{gl}_{2n}$ lies in
$\liea{sp}_{2n}$ if and only if $(Jx,z)=-(x,Jz)$ or, filling in the
expression for $J$, if
\[ \langle x, \eta \rangle (y,z)=-(x,y) \langle z,\eta \rangle \]
for all $x,z \in V$. By skewness of the form, this is clearly the case
if $\eta=\phi(y)$. Conversely, for $z \in V$ fixed such that $(y,z)
\neq 0$, the equation above shows that $\eta=c\phi(y)$ for some $c \neq
0$; hence if $d$ is a square root of $c$, then $J=y \otimes \eta= dy
\otimes \phi(dy)$. This shows that the second set of the lemma is equal
to the third, and that they contain $C$. Finally, the transitivity of
$\lieg{Sp}_{2n}$ on $V \setminus \{0\}$ implies the transitivity of
$\lieg{Sp}_{2n}$ on the second set of the lemma, and this concludes
the proof.
\end{proof}

\begin{proof}[Proof of Theorem \ref{thm:slsp} for $\liea{sp}_{2n}$]
Let $A \in \liea{sp}_{2n}$ be non-zero. We show that there exists an $x
\in V$ such that $A-(Ax \otimes \phi(Ax))$ has rank $\rk(A)-1$. By Lemma
\ref{lm:Rank1} this is the case if and only if $\ker \phi(Ax) \supseteq
\ker A$ and $\langle x, \phi(Ax) \rangle=1$. The first condition holds
for all $x \in V$, as $Az=0$ implies
\[ \langle z,\phi(Ax) \rangle=(z,Ax)=-(Az,x)=0 \]
by virtue of $A^T=-A$. Hence, we are left to show that there exists
an $x \in V$ for which $\langle x,\phi(Ax) \rangle=(x, Ax) \neq 0$;
rescaling $x$ will then make this scalar $1$. Suppose, on the contrary,
that $(x, Ax)=0$ for all $x$. Then we have for all $x,y \in V$:
\[ 0=(x+y,A(x+y))=(x,Ax)+(y,Ay)+(x,Ay)+(y,Ax)=2(x,Ay), \]
so that $A=0$, which contradicts our assumption that $\rk(A)$ be
greater than $1$.  

By induction, this shows that for $k=1,\ldots,2n$ every element
of $\liea{sp}_{2n}$ of rank $k$ lies in $kC$. As in the case of
$\liea{sl}_n$, the inclusions $kC \subseteq (k+1)C$ and $0 \in 2C$,
together with the fact that $kC$ cannot contain elements of rank higher
than $k$, conclude the proof.
\end{proof}

\section{The main course: $\liea{o}_n$} \label{sec:o}

Now we come to the more intricate part of this paper: the secant
varieties of the minimal orbit $C$ of $\lieg{SO}_n$ on its Lie algebra
$\liea{o}_n$. Unlike in the cases of $\liea{sl}_n$ and $\liea{sp}_{2n}$,
the sums $kC$ are in general not closed, and only their closures are
described explicitly here. The approach, though, is the same as for
$\liea{sl}_n$ and $\liea{sp}_{2n}$: we try to decrease the rank of a
given element of $\liea{o}_n$ by subtracting an appropriate element
of $C$. How this rank reduction works for $\liea{o}_n$, and why it
comes short of characterising the sets $kC$ completely, is explained in
Subsection \ref{ssec:C}. Subsections \ref{ssec:2C}, \ref{ssec:3C}, and
\ref{ssec:kC} are devoted to determining $\overline{2C}$, $\overline{3C}$,
and $\overline{kC}$ for $k \geq 4$, respectively. In Subsection
\ref{ssec:2C} we find that the complement of $2C$ in $\overline{2C}$
is a single nilpotent orbit, which discovery motivates the discussion
of nilpotent orbits in Subsection \ref{ssec:Nilp}.

\subsection{The minimal orbit and rank reduction} \label{ssec:C}

We retain the notation $\phi$ and $A^T$ from Section \ref{sec:slsp}; only
now they are defined with respect to the non-degenerate {\em symmetric}
bilinear on $V=K^n$ defining the Lie algebra $\liea{o}_n$. Recall that,
for any $A \in \liea{o}_n$ and $\lambda \in K$, the numbers $\lambda$
and $-\lambda$ have the same (geometric and algebraic) multiplicity among
the eigenvalues of $A$; moreover, $\rk(A)$ is even. The following
lemma implies that every $A \in \liea{o}_n$ of rank $2k$ is the sum of
$k$ rank two elements of $\liea{o}_n$.

\begin{lm} \label{lm:RankRed2}
Let $A \in \liea{o}_n$, and let $J \in \liea{gl}_n$ of rank one be such
that $\rk(A-J)=\rk(A)-1$. Then $\rk(A-(J-J^T))=\rk(A)-2$.
\end{lm}

The proof of this lemma uses the useful identities $(y \otimes
\eta)^T=\phi^{-1} \eta \otimes \phi y\ (y \in V, \eta \in V^*)$
and $\phi A^T=A^* \phi\ (A \in \liea{gl}_n)$, whose proofs are
straightforward.

\begin{proof}
By Lemma \ref{lm:Rank1}, there exist $x \in V$ and $\xi \in V^*$ such
that $J=Ax \otimes A^*\xi$ and $\la Ax,\xi \ra=1$; note that then
$\ker(A-J)=\ker A \oplus Kx$.  We have
\[ J^T 
= \phi^{-1} A^* \xi \otimes \phi Ax 
= A^T \phi^{-1} \xi \otimes (A^T)^* \phi x
= A\phi^{-1} \xi \otimes A^* \phi x,
\]
where the third step is justified by $A^T=-A$. In particular, we find
that $\ker J^T$, too, contains $\ker A$, so that $\ker A \subseteq \ker
(A-J+J^T)$. Moreover, we have
\[
        (A-J+J^T)x=J^T x=\la x, A^* \phi x \ra A \phi^{-1} \xi
        =(Ax,x) A \phi^{-1} \xi=0,
\]
while $Ax \neq 0$. (In the last step we used $(Ax,x)=(x,A^T
x)=-(x,Ax)=-(Ax,x)$.) Hence, $\rk(A-J+J^T)$ is strictly smaller than
$\rk(A)$; but as $A-J+J^T$ is skew symmetric, its ranks is even, hence
equal to $\rk(A)-2$.
\end{proof}

If $C$ would contain all elements of $\liea{o}_n$ of rank $2$, then we
would have $mC=\liea{o}_n$ by this lemma. However, $C$ is smaller; to
characterise it we first describe the rank-two-elements of $\liea{o}_n$.

\begin{prop} \label{prop:Rank2}
For any $2$-dimensional subspace $W=\langle y_1, y_2 \rangle_K$ of $V$,
the space $\{A \in \liea{o}_n \mid \im A \subseteq W \}$ is one-dimensional and
spanned by $y_1 \otimes \phi(y_2) - y_2 \otimes \phi(y_1)$.
\end{prop}

The proof of this proposition uses another easy observation; namely, that for
any $A \in \liea{o}_n$ the kernel of $A$ is the orthogonal complement
of $\im A$ with respect to $(.,.)$; we denote this orthogonal complement
by $(\im A)^\bot$.

\begin{proof}
Let $A \in \liea{o}_n \setminus \{0\}$ have image contained in, and hence
equal to, $W$; and let $x_1 \in V$ be such that $Ax_1=y_1$. Then we have
$(x_1,y_1)=(x_1,Ax_1)=0$ by the skewness of $A$, so that $(x_1,y_2)=0$
would imply $x_1 \in (\im A)^\bot=\ker A$, a contradiction, hence we may
set $\alpha:=1/(x_1,y_2)$. Furthermore, $y_2^\bot \supseteq \ker A$,
so that $J:=\alpha y_1 \otimes \phi (y_2) \in \liea{gl}_n$ satisfies
the condition of Lemma \ref{lm:RankRed2}. Then that lemma implies
$A= \alpha (y_1 \otimes \phi(y_2) - y_2 \otimes \phi(y_1))$, as
claimed.
\end{proof}

Proposition \ref{prop:Rank2} has the following interesting
consequence.

\begin{cor} \label{cor:Rank2}
For each $k \in \{0,1,2\}$, the group $\lieg{SO}_n$ acts transitively
on the set $O_k:=\PP\{A \in \liea{o}_n \mid \rk(A)=2 \text{ and } 
	(.,.)|_{\im A} \text{ has rank } k \} \subseteq \PP \liea{g}$.
\end{cor}

\begin{proof}
It is not hard to see that $\lieg{SO}_n$ acts transitively on the
$2$-dimensional subspaces of $V$ on which $(.,.)$ has rank $k$, and now
the proposition can be applied.
\end{proof}

The following corollary identifies $\PP C$ with $O_0$.

\begin{cor} \label{cor:Con}
The set $C$ consists of all $A \in \liea{o}_2$ with $\rk(A) \leq 2$ and
$\im A$ isotropic with respect to $(.,.)$. The latter condition is
equivalent, for $A \in \liea{o}_2$, to $A^2=0$.
\end{cor}

\begin{proof}
By Corollary \ref{cor:Rank2}, it suffices to check that the highest root
vector of $\liea{o}_n$ with respect to some choice of Borel and Cartan
subalgebras has the stated properties, which is straightforward. As
for the second statement: the radical of $(.,.)|_{\im A}$ on $\im A$ is
exactly $\ker A \cap \im A$, hence all of $\im A$ if and only if $A^2=0$.
\end{proof}

We are now ready to state and prove our main rank reduction argument
in the orthogonal case.

\begin{prop} \label{prop:RankRed}
Let $A \in \liea{o}_n$ be of rank $\geq 4$. Then there exists a $J \in C$
such that $\rk(A-J)=\rk(A)-2$.
\end{prop}

\begin{proof}
On $\im A$ we have two bilinear forms: the restriction of $(.,.)$, and
a second form $(.|.)$ defined by $(Ax_1|Ax_2)=(x_1,Ax_2)$; we continue
to use $\bot$ only for `perpendicular with respect to $(.,.)$'.
The second form is well-defined as $\ker A \bot \im A$ and
skew-symmetric because
\[ (Ax_2|Ax_1)=(x_2,Ax_1)=(A^T x_2,x_1)=-(Ax_2,x_1)=-(x_1,Ax_2)=-(Ax_1|Ax_2). 
\]
Moreover, $(.|.)$ is non-degenerate, as $(Ax_1|Ax_2)=0$ for all $x_1$ implies
$x_1 \bot Ax_2$ for all $x_1$, i.e., $Ax_2=0$.  We may now apply Lemma
\ref{lm:BilForms} below to find a $2$-dimensional subspace $U$ 
of $\im A$ that is isotropic with respect to $(.,.)$ but not with
respect to $(.|.)$. Choose a basis $y_1,y_2$ of $U$ such that
$(y_1|y_2)=1$, and set $J:=y_1 \otimes \phi(y_2)-y_2 \otimes
\phi(y_1)$. Then $\im J$ is two-dimensional and isotropic with
respect to $(.,.)$, so $J$ lies in $C$ by Corollary \ref{cor:Con}.
Furthermore, $\ker J=U^\bot$ contains $\ker A=(\im A)^\bot$, and if
$x_1 \in A^{-1}y_1$, then 
\[ (A-J)x_1=y_1-(x_1,y_2)y_1+(x_1,y_1)y_2=
	y_1-(y_1|y_2)y_1+(y_1|y_1)y_2=0,
\]
so that $\rk(A-J)$ is strictly smaller than $\rk(A)$; we conclude that
$J$ has the required properties.
\end{proof}

The proof above uses the following observation on bilinear forms.

\begin{lm} \label{lm:BilForms}
Let $W$ be a $K$-vector space of finite dimension $\geq 4$ equipped with a
(possibly degenerate) symmetric bilinear form $B_1$ and a non-degenerate
skew-symmetric bilinear form $B_2$. Then there exists a $2$-dimensional
subspace of $W$ that is isotropic with respect to $B_1$ but not with
respect to $B_2$.
\end{lm}

The following proof, which is considerably shorter than our original
proof, was suggested by Jochen Kuttler.

\begin{proof}
Suppose, on the contrary, that all $2$-dimensional $B_1$-isotropic
subspaces of $W$ are $B_2$-isotropic, and note that then
all $B_1$-isotropic subspaces of {\em any} dimension are
$B_2$-isotropic. We may choose a basis $e_1,\ldots,e_d$ of $W$ such
that $B_1(x,y)=\sum_{j=1}^l x_j y_j$, where $l \leq d$ is the rank of
$B_1$. If $l=0,1$, or $2$, then the subspace of codimension $1$ defined
by the equation $x_1=0,x_1=0$, or $x_2=ix_1$, respectively, is isotropic
with respect to $B_1$, and hence with respect to $B_2$. On the other hand,
any $B_2$-isotropic subspace of $W$ has dimension at most $\dim(W)/2$,
so that $\dim(W)-1 \leq \dim(W)/2$, a contradiction to $\dim(W) \geq 4$.

Hence $l$ is at least $3$. Now consider the quadric $Q_1:=\{x \in W \mid
B_1(x,x)=0\}$.  It is easy to find linearly independent vectors $w_1,
w_2, w_3$ on $Q_1$ such that $B_1(w_j,w_k)\neq 0$ for all distinct
$j,k \in \{1,2,3\}$---for example, $w_1=e_1+ie_2,w_2=e_1+ie_3,$ and
$w_3=e_2+ie_3$. Moreover, $Q_1$ spans $W$ and we may find $w_4,\ldots,w_d
\in Q_1$ such that $w_1,w_2,\ldots,w_d$ is a basis of $W$; we write
$w_j^\bot$ for $\{w \in W \mid B_1(w_j,w)=0\}$. For each $j$ and any $w
\in Q_1 \cap w_j^\bot$ the space $Kw_j + Kw$ is $B_1$-isotropic, so that
$B_2(w_j,w)=0$ by assumption. As, moreover, the restriction of $B_1$
to $w_j^\bot$ has rank at least two, we find that $Q_1 \cap w_j^\bot$
spans $w_j^\bot$, so that $B_2(w_j,w)=0$ for any $w \in w_j^\bot$. In
other words, the linear function $B_2(w_j,.)$ is equal to $c_j B_1(w_j,.)$
for some $c_j \in K$, so that if $A_1,A_2$ are the matrices of $B_1,B_2$
with respect to $w_1,\ldots,w_d$, then
\[ A_2=\diag(c_1,\ldots,c_d) A_1. \]
As $A_2$ is skew (with respect to transposition in the main diagonal)
and $A_1$ is symmetric, we find that $c_j a_{jk}=-c_k a_{jk}$ for all
$j,k=1,\ldots,d$. By construction $a_{jk} \neq 0$ for distinct $j,k
\in \{1,2,3\}$, and we find that $c_j=-c_k$ for all such $j,k$. This
readily implies that $c_1=c_2=c_3=0$, so that $A_2$ is singular; but
this contradicts the non-degeneracy of $B_2$.
\end{proof}

We can now prove Theorem \ref{thm:o2}; from the proof it will become
clear why the rank reduction of Proposition \ref{prop:RankRed} does
not suffice to characterise the secant varieties of $C$ completely.

\begin{proof}[Proof of Theorem \ref{thm:o2}.]
By Proposition \ref{prop:RankRed} and induction, it suffices to prove
that every element of $\liea{o}_n$ having rank $2$ lies in $4C$. By
Corollary \ref{cor:Rank2} (and the fact that $4C$ is, of course, a
cone) it suffices to prove this for particular representatives of the
projective orbits $O_k\ (k=0,1,2)$ mentioned in that corollary. For
$k=0$ we have $O_0=\PP C$, so there is nothing to prove. For $k=1,2$
let $y_1,y_2,y_3,y_4$ be linearly independent isotropic vectors in $V$
satisfying $(y_1,y_3)=(y_2,y_4)=1$ and $\la y_1, y_3 \ra_K \bot \la y_2,
y_4 \ra_K$ (such vectors exist). Then a representative of $O_1$ is
\[ (y_1+y_3) \otimes \phi(y_2)-y_2 \otimes \phi(y_1+y_3), \]
which can be written as
\[ (y_1 \otimes \phi(y_2) - y_2 \otimes \phi(y_1))
        +(y_3 \otimes \phi(y_2) - y_2 \otimes \phi(y_3)) \in 2C.
\]
Similarly, a representative of $O_2$ is
\[ (y_1+y_3) \otimes \phi(y_2+y_4)-(y_2+y_4) \otimes \phi(y_1+y_3),
\]
which equals
\begin{align*} &(y_1 \otimes \phi(y_2) - y_2 \otimes \phi(y_1))
        +(y_1 \otimes \phi(y_4) - y_4 \otimes \phi(y_1))\\
        +&(y_3 \otimes \phi(y_2) - y_2 \otimes \phi(y_3))
        +(y_3 \otimes \phi(y_4) - y_4 \otimes \phi(y_3)) \in 4C.
\end{align*}
\end{proof}

One may think, now, that a representative of $O_2$ could already lie in
$kC$ for $k=2$ or $3$---but this is not the case. Indeed, as we shall
see in Subsection \ref{ssec:3C}, such a representative does not even
lie in $\overline{3C}$. This serves to show that the secant varieties
of the minimal orbit in $\liea{o}_n$ are considerably more complicated
than those of the minimal orbits in $\liea{sl}_n$ and
$\liea{sp}_{2n}$.

\subsection{The first secant variety} \label{ssec:2C}

The first secant variety $\PP(\overline{2C}\setminus \{0\})$ of
the minimal orbit in {\em any} simple Lie algebra is described in
\cite{Kaji2000} as the union of a single (projective) semisimple orbit and
several nilpotent orbits. We reprove this statement here for $\liea{o}_n$;
first, because our method is different from that of Kaji {\em et al}
and also applies to the second secant variety, and second, because we
want to determine the complement $\overline{2C}\setminus 2C$ explicitly.

Before stating our characterisation of $\overline{2C}$, we recall that
the closed subvariety
\[ R_k := \{ A \in \liea{o}_n \mid \rk A \leq 2k \} \]
is irreducible for all $k$---this follows, for instance, from 
\cite[Lemma 4.2.4(3)]{Goodman98}---and we recall from the
introduction the notation $S_2$ for the set of rank $4$ semisimple
elements having non-zero eigenvalues $a,a,-a,-a$.

\begin{prop} \label{prop:2Cbar}
The affine variety 
\[ M:=\{A \in \liea{o}_n \mid \rk(A) \leq 4 \text{ and } A^3=\lambda A
	\text{ for some } \lambda \in K\} \]
has two irreducible components, namely $R_1$ and $\overline{2C}$.
Furthermore, $\overline{2C}$ is equal to $\overline{S_2}$.
\end{prop}

Recall, for the proof of this proposition, the notation $\mathcal{O}[\bd]$
for the nilpotent $\lieg{O}_n$-orbit on $\liea{o}_n$ corresponding to
the partition $\bd$ of $n$, where the even entries of $\bd$ are supposed
to have even multiplicities. We work with $\lieg{O}_n$ here, rather than
with the adjoint group $\lieg{SO}_n$, not to have to distinguish between
the two $\lieg{SO}_n$-orbits corresponding to {\em very even} partitions
\cite{Collingwood93,Kraft82}. Indeed, as both groups have the same minimal
orbit $C=\mathcal{O}[2,2,1^{n-4}]$, this subtlety is immaterial to us.

We will not be able to avoid, in what follows, some explicit matrix
computations. In these computations we always take for $(.,.)$ the
symmetric form given by $(x,y)=\sum_{i=1}^n x_i y_{n+1-i}$ with respect
to the standard basis of $V=K^n$. The elements of $\liea{o}_n$ are then
skew symmetric about the {\em skew} diagonal running from position $(1,n)$
to position $(n,1)$.

\begin{proof}[Proof of Proposition \ref{prop:2Cbar}.]
For $J_1,J_2 \in C$ we have
\[ (J_1+J_2)^3=J_1J_2J_1+J_2J_1J_2, \]
where we use that $J_i^2=0$ for $i=1,2$ (Corollary \ref{cor:Con}).
The map $J_1J_2J_1$ is skew-symmetric and its image is contained in
$\im J_1$, hence by Proposition \ref{prop:Rank2} $J_1J_2J_1=c_1 J_1$ for
some $c_1 \in K$. Similarly, $J_2J_1J_2=c_2 J_2$ for some $c_2 \in K$. If
$J_1J_2=0$, then $c_1=c_2=0$ and $J_1+J_2 \in M$ (with $\lambda=0$).
Otherwise, let $x \in V$ be such that $J_1J_2x \neq 0$. Then
\[ c_2 J_1 J_2 x = J_1(J_2 J_1 J_2)x = (J_1 J_2 J_1) J_2 x= c_1 J_1 J_2 x, \] 
so that $c_1=c_2$. This shows that $J_1+J_2$ lies in $M$ (with
$\lambda=c_1$), so that $\overline{2C} \subseteq M$. The inclusion $R_1
\subseteq M$ is immediate: an element $A$ of $R_1$ is either semisimple
with non-zero eigenvalues $a,-a$, so that $A \in M$ (with $\lambda=a^2$),
or it is nilpotent of nilpotence degree at most $3$, and then $A$ also
lies in $M$ (with $\lambda=0$).

Conversely, let $A$ be in $M$ and let $\lambda \in K$ be such that
$A^3 = \lambda A$. If $\lambda=0$, then $A^3=0$, which together with
the condition that $\rk(A)$ be at most $4$ shows that $A$ lies in a
nilpotent orbit corresponding to one of the partitions $[3,3,1^{n-6}]$,
$[3,2,2,1^{n-7}]$, $[3,1^{n-3}]$, $[2,2,2,2,1^{n-8}]$, $[2,2,1^{n-4}]$,
or $[1^n]$. The first among these is greater than all of the other five
in the usual order on partitions \cite{Collingwood93}, so that the
corresponding orbit closure contains the other five nilpotent orbits.

Suppose, on the other hand, that $\lambda \neq 0$, and let $a$ be a
square root of $\lambda$. Then $A$ is a zero of the square-free polynomial
$t(t-a)(t+a)$, hence semisimple. There are three possibilities: either
$A=0$, or $A \in R_1$ with non-zero eigenvalues $\pm a$, or $A \in S_2$
with eigenvalues $a,a,-a,-a$. Together with the above discussion
of the nilpotent orbits in $M$ this implies
$M= R_1 \cup S_2 \cup \overline{\mathcal{O}}[3,3,1^{n-6}]$.
We shall show that the last two terms are contained in $\overline{2C}$,
so that
\[ M= R_1 \cup \overline{2C}; \]
as $R_1$ and $\overline{2C}$ are both irreducible and neither of these
sets is contained in the other, this implies the first statement
of the proposition. Moreover, the above shows that the complement
of $S_2$ in $\overline{2C}$ equals $(R_1 \cap \overline{2C}) \cup
\overline{\mathcal{O}}[3,3,1^{n-6}]$, so that $S_2$ is open, and hence
dense, in $\overline{2C}$---which proves the second statement of the
proposition.

Suppose, therefore, that $A$ lies in $\mathcal{O}[3,3,1^{n-6}]$. Then $A$
is conjugate to an $n \times n$-matrix that has a $6 \times 6$ block
\begin{small}
\[ 
\begin{bmatrix}
0&1&&&&\\
&0&1&&&\\
&&0&&&\\
&&&0&-1&\\
&&&&0&-1\\
&&&&&0
\end{bmatrix}
=
\begin{bmatrix}
0&1&&&&\\
&0&&&&\\
&&0&&&\\
&&&0&&\\
&&&&0&-1\\
&&&&&0
\end{bmatrix}
+
\begin{bmatrix}
0&&&&&\\
&0&1&&&\\
&&0&&&\\
&&&0&-1&\\
&&&&0&\\
&&&&&0
\end{bmatrix} 
\]
\end{small}
\!\!\!\!\!
in the middle, and zeroes elsewhere (off-diagonal zeroes are omitted). Now
both matrices on the right-hand side of the equality lie in $C$: they
have rank $2$ and isotropic images.  Therefore, $A$ lies in $2C$ and
$\overline{\mathcal{O}}[3,3,1^{n-6}] \subseteq 2C$.

Next assume that $A$ has rank $4$ and is semisimple with non-zero
eigenvalues $a,a,-a,-a$. Then $A$ is conjugate to a matrix with
a $4 \times 4$-block
\begin{small}
\[ 
a
\begin{bmatrix}
1&&&\\
&1&&\\
&&-1&\\
&&&-1
\end{bmatrix}
=
\frac{a}{2}
\begin{bmatrix}
1&&1&\\
&1&&-1\\
-1&&-1&\\
&1&&-1
\end{bmatrix}
+\frac{a}{2}
\begin{bmatrix}
1&&-1&\\
&1&&1\\
1&&-1&\\
&-1&&-1
\end{bmatrix}
\]
\end{small}
\!\!\!\!
in the middle and zeroes elsewhere. One readily verifies that the
two terms on the left-hand side lie in $C$, so that $A \in 2C$. As
explained above, this concludes the proof of the proposition.
\end{proof}

Now that we have identified $\overline{2C}$ as the irreducible component
$\overline{S_2}$ of $M$---and hence proved the first part of Theorem
\ref{thm:o}---we investigate the set $2C$ itself. It is easy to verify,
like we did in the proof above for $\mathcal{O}[3,3,1^{n-6}]$, that the
nilpotent orbits $\mathcal{O}[3,1^{n-3}]$, $\mathcal{O}[2,2,2,2,1^{n-8}]$,
$\mathcal{O}[2,2,1^{n-4}]$, and $\mathcal{O}[1^n]$ in $\overline{S_2}$
all lie in $2C$ (in fact, this follows from the computations in Subsection
\ref{ssec:Nilp} below), as does $S_2$ by the explicit computation in the
proof above. Thus we find that $\overline{2C} \setminus 2C$ is contained
in $\mathcal{O}[3,2,2,1^{n-7}]$. Note that an element $A$ from this
nilpotent orbit has $\rk(A)=4$ and $\rk(A^2)=1$. The following lemma
shows that an element of $2C$ cannot have this property, thus proving
Proposition \ref{prop:2C}.

\begin{lm} \label{lm:rk-even}
If $A \in 2C$ has rank $4$, then $\rk(A^2)$ is even.
\end{lm}

\begin{proof}
By Corollary \ref{cor:Con} we may write 
\[ A=y_1 \otimes \phi(y_2) - y_2 \otimes \phi(y_1) +
	y_3 \otimes \phi(y_4) - y_4 \otimes \phi(y_4), \]
where $\langle y_1,y_2 \rangle_K$ and $\langle y_3,y_4 \rangle_K$
are isotropic. By the condition that $\rk(A)$ be $4$, the vectors
$y_1,y_2,y_3,y_4$ are a basis of $\im A$. The matrix of
$(.,.)|_{\im A}$ with respect to this basis is of the form 
\[ \begin{bmatrix} 0 & M \\ M^t &0 \end{bmatrix}, \]
where $M^t$ is the (ordinary) transpose of the $2\times 2$-matrix $M$.
We conclude that $\rk A^2=\rk (.,.)|_{\im A}=2\rk M$.
\end{proof}

\subsection{The second secant variety} \label{ssec:3C}

To characterise $\overline{3C}$ we proceed as in the first part of
the proof of Proposition \ref{prop:2Cbar}: we take three arbitrary
elements $J_1,J_2,J_3$ of $C$, and use the relations provided
by Proposition \ref{prop:Rank2} to find a polynomial annihilating
$J_1+J_2+J_3$. Conversely, we show that semisimple elements having a
characteristic polynomial of that form do indeed lie in $3C$.

\begin{prop} \label{prop:3C}
Any element of $\overline{3C}$ is annihilated by a polynomial of the form
\[ t (t-a)(t-b)(t-c)(t+a)(t+b)(t+c) \]
for some $a,b,c \in K$ with $a+b+c=0$.
\end{prop} 

\begin{proof}
The set of matrices in $\liea{o}_n$ that are annihilated by such a
polynomial, is closed, so that it suffices to prove the proposition for
elements of $3C$. Let therefore $J_1,J_2,J_3$ be elements of $C$. From
the proof of Proposition \ref{prop:2Cbar} we know that there exist constants
$c_{ik}=c_{ki}$ such that
\[ J_i J_k J_i = c_{ik} J_i \text{ for all } i,k \in \{1,2,3\},\ i \neq k. \]
Similarly, there exists a $c \in K$ with
\[ J_1 (J_2J_3J_1J_2J_3+J_3J_2J_1J_3J_2) J_1 = c J_1; \]
indeed, the matrix between brackets is an element of $\liea{o}_n$, so
that the matrix on the left-hand side lies in $\liea{o}_n$. Furthermore,
its image is contained in $\im J_1$, whence the existence of such a
$c$ follows from Proposition \ref{prop:Rank2}. In fact, one can show
that the same $c$ satisfies the above relation with $1,2,3$ permuted
cyclically. Using these relations, a straightforward calculation
shows that
\[ t p(t) \text{ with }
        p(t):=(t^3-(c_{12}+c_{13}+c_{23})t)^2-c-2c_{12}c_{13}c_{23}
\]
annihilates $J_1+J_2+J_3$. Now we need only check that $p$ has the desired
form. To this end, let $\mu$ be a square root of $c+2c_{12}c_{13}c_{23}$,
so that $p$ factorises into
\[ p(t)=(t^3-(c_{12}+c_{13}+c_{23})t + \mu) (t^3-(c_{12}+c_{13}+c_{23})t - \mu).
\]
The first of these factors lacks a term with monomial $t^2$; hence, the
sum of its zeroes $a,b,c$ is $0$. The second factor has zeroes $-a,-b,-c$,
and this concludes the proof of the proposition.
\end{proof}

\begin{re}
The polynomial $tp(t)$ appearing in the proof above was found as follows:
consider the free associative algebra $F$ (with one) over the ground field
$K(c_{12},c_{13},c_{23},c)$ with generators $J_1,J_2,J_3$ and let $I$ be
the ideal generated by the relations appearing in the proof above.  Then a
(non-commutative) Gr\"obner basis computation of $I$ shows that $F/I$
has dimension $37$, and the polynomial $tp(t)$ is the minimal polynomial
of $J_1+J_2+J_3$ in this quotient. For this computation we used the {\tt
GAP}-package {\tt GBNP} written by Cohen and Gijsbers \cite{Cohen2003, GAP4}
(with concrete values for the $c_{ij}$ and $c$), together with some {\em
ad hoc} programming of our own in {\tt Mathematica}.
\end{re}

A partial converse to the proposition above is the following lemma,
in whose proof we compute with respect to the fixed bilinear form of
Subsection \ref{ssec:2C}.

\begin{lm} \label{lm:3C}
For all $a,b \in K$, any semisimple element of $so_n$ whose eigenvalues
(with multiplicities) are $0$ ($n-6$ times) and $a,b,-a-b,-a,-b,a+b$
for some $a,b \in K$, lies in $3C$.
\end{lm}

\begin{proof}
Let $A$ be such an element; we may suppose that $A$ is non-zero. Then
the numbers $a,b,-a-b$ are not all equal, and by permuting them we may
assume that $a \neq b$. Now $A$ is conjugate to an $n \times n$-matrix
having zeroes everywhere except for a $6 \times 6$-block in the
middle, which is of the form
\begin{small}
\begin{align*}
&\begin{bmatrix}
a & 0 & 0 & 0 & 0 & 0\\
0 & b & 0 & 0 & 0 & 0\\
0 & 0 &a+b& 0 & 0 & 0\\
0 & 0 & 0&-a-b& 0 & 0\\
0 & 0 & 0 & 0 &-b & 0\\
0 & 0 & 0 & 0 & 0 &-a\\
\end{bmatrix}
=\frac{1}{b-a}
\begin{bmatrix}
ab &ab &0  &0  &0  &0\\
-ab&-ab &0  &0  &0  &0\\
0  &0  &0  &0  &0  &0\\
0  &0  &0  &0  &0  &0\\
0  &0  &0  &0 &ab &-ab\\
0  &0  &0  &0 &ab &-ab\\
\end{bmatrix}\\
&+\frac{1}{b-a}
\begin{bmatrix}
-a^2&-ab &0  &0  &0  &0\\
ab &b^2 &0  &0  &0  &0\\
0  &0  &(b+a)(b-a)  &0  &0  &0\\
0  &0  &0 &-(b+a)(b-a)  &0  &0\\
0  &0  &0  &0&-b^2 &ab\\
0  &0  &0  &0&-ab &a^2\\
\end{bmatrix}.
\end{align*}
\end{small}
\!\!\!\!\!\!
The first term on the right-hand side lies in $C$, as it has rank
$2$ and its image is isotropic. We claim that the second term on the
right-hand side, which we denote by $B$, lies in $2C$. To see this, note
that $(b,-a,0,0,0,0)^t$, $(0,0,1,0,0,0)^t$, and $(a,-b,0,0,0,0)^t$ are
eigenvectors of $B$ with eigenvalues $0,a+b,$ and $a+b$, respectively. It
follows that $B$ has rank $4$ and that $-a-b$, as well, has multiplicity
two among the eigenvalues of $B$; hence $B$ lies in $2C$ by
Propositions \ref{prop:2C} and \ref{prop:2Cbar}, and $A$ lies in $3C$.
\end{proof}

To finish our characterisation of $\overline{3C}$, recall that $R_k$
is the set of all elements of $\liea{o}_n$ having rank $\leq 2k$. We
now need an argument why $S_3$, the set of all semisimple elements in
$R_3$ having $6$ {\em distinct} eigenvalues $a,b,c,-a,-b,-c$ such that
$a+b+c=0$, is dense in $\overline{3C}$. The following lemma will provide
such an argument.

\begin{lm} \label{lm:Tk}
For any $k$, the subset $T_k$ of $R_k$ consisting of all elements having
$2k$ distinct non-zero eigenvalues, is open in $R_k$.
\end{lm}

\begin{proof}
An element of $R_k$ lies in $R_k \setminus T_k$ if and only if it has
a characteristic polynomial of the form
\[ t^{n-2k} (t^2-a_1)^2(t^2-a_2)(t^2-a_3)\ldots(t^2-a_{k-1}) \]
for some $a_1,\ldots,a_{k-1} \in K$. The map $K^{k-1} \rightarrow K^{n+1}$
sending $(a_1,\ldots,a_{k-1})$ to the coefficients of the monomials $t^i$
in the polynomial above has a closed image $Y$, and $R_k \setminus T_k$
is the inverse image of $Y$ under the polynomial map sending a matrix
to the coefficients of its characteristic polynomial.
\end{proof}

From Proposition \ref{prop:3C} and Lemma \ref{lm:3C} we have
\[ S_3=T_3 \cap \overline{3C}. \]
By the lemma above, this set is open, and hence dense, in $\overline{3C}$,
so that $\overline{3C}=\overline{S_3}$ as claimed in Theorem \ref{thm:o}.

\subsection{Higher secant varieties} \label{ssec:kC}

After reading the discussion of $\overline{2C}$ and $\overline{3C}$,
one could think that to describe the sets $\overline{kC}$ for
$k \geq 4$, we must consider the quotient of the free algebra
generated by $J_1,\ldots,J_k$ by the ideal generated by all
relations that can be inferred from Proposition \ref{prop:Rank2},
i.e., those of the form $J_i J_k J_i=c_{ik} J_i$ appearing in
the proof of Proposition \ref{prop:2Cbar}, those reflecting that
$J_i(J_kJ_lJ_iJ_kJ_l+J_lJ_kJ_iJ_lJ_kJ_i)J_i$ is a scalar multiple of $J_i$
(where the scalar does not change if we permute $i,k,l$ cyclically; this
relation appears in the proof of Proposition \ref{prop:3C}), and similar
relations, such as: $J_i(J_kJ_lJ_m+J_mJ_lJ_k)J_i$ is a scalar multiple
of $J_i$. While this quotient algebra may be interesting in itself---is
it always finite-dimensional? Gr\"obner Basis computations seem to end
in an endless loop already for $k=4$---it does, surprisingly enough,
not play an important role in determining the higher secant varieties
of $C$. The following proposition explains why.

\begin{prop}
The set $4C$ contains a dense subset of $R_4$.
\end{prop}

Together with the obvious inclusion $4C \subseteq R_4$, this proposition
implies $\overline{4C}=R_4$. Using Proposition \ref{prop:RankRed}
we then find $\overline{kC}=R_{2k}$ for all $k \geq 4$, as claimed in
Theorem \ref{thm:o}.

\begin{proof}
Let $a_1,a_2,a_3,a_4$ be variables. It suffices to prove that the
diagonal matrix
\[ A=\diag\{a_1,\ldots,a_4,-a_4,\ldots,-a_1\} \]
lies in $4C(K(a_1,\ldots,a_4))$, i.e., $4$ times the minimal orbit in
$\liea{o}_8$ with coordinates in $K(a_1,\ldots,a_4)$. Indeed, if this
is the case, then a generic semisimple element of $R_4$ lies in $4C$,
and these elements are dense in $R_4$ by Lemma \ref{lm:Tk}. Define
the expressions
\begin{align*}
r_1&:=0,
&s_1&:=1,\\
r_2&:=a_3(a_1^2-(a_2+a_3+a_4)^2),
&s_2&:=4(a_2+a_3)(a_3+a_4),\\
r_3&:=-a_4(a_1^2-(-a_2+a_3+a_4)^2),
&s_3&:=4(a_3+a_4)(-a_2+a_4),\\
r_4&:=-a_2(a_1^2-(-a_2-a_3+a_4)^2), \text{ and}
&s_4&:=4(-a_2+a_4)(-a_2-a_3);
\end{align*}
and note that the transformation $a_2 \mapsto a_3 \mapsto a_4 \mapsto
-a_2$ cyclically permutes $s_2,s_3,s_4$, and does the same with
$r_2,r_3,r_4$ up to a change of sign. Now set 
\begin{align*}
y_1&:=(0,r_4,0,r_3,0,r_2,0,r_1)^t,&
y_2&:=(\frac{1}{s_1},0,\frac{1}{s_2},0,\frac{1}{s_3},0,\frac{1}{s_4},0)^t,
	\text{ and}\\
J&:=y_1 (y_2^t F) - y_2 (y_1^t F),
\end{align*}
where $F=(\delta_{i+j,9})_{ij}$ is the $8 \times 8$-matrix
representing the form $(.,.)$. By construction $J$ lies in
$\liea{o}_n(K(a_1,\ldots,a_4))$ and has rank $2$. It is easy to see that
$y_1^t F y_1=y_2^t F y_2=0$, and a straightforward computation shows that
$y_1^t F y_2 = \sum_{i=1}^4 \frac{r_i}{s_i}$ is zero, as well. This shows
that $\im J$ is isotropic, hence $J$ lies in $C(K(a_1,\ldots,a_4))$. A
direct computation (preferably by a computer algebra system; we used
{\tt Mathematica}) shows that $A-J$ is semisimple with eigenvalues
\[ 0,0,\mp a_1,\pm\frac{1}{2}(a_1-a_2+a_3-a_4),
	\pm\frac{1}{2}(a_1+a_2-a_3+a_4); \]
if we take the upper one of the two signs in each of the last
three eigenvalues, then they add up to zero, so that $A-J \in
3C(K(a_1,\ldots,a_4))$ by Lemma \ref{lm:3C}. We conclude that $A$ lies
in $4C(K(a_1,\ldots,a_4))$ as claimed.
\end{proof}

\begin{re} 
By studying the computations needed for the proof above, it should be
straightforward to prove that {\em any} semisimple element of $R_4$
lies in $4C$. Furthermore, the computation proving Lemma \ref{lm:3C} is
easily modified to a proof that for any semisimple $A \in \liea{o}_n$
of rank $2k,\ k \geq 2$, there exists an element $J \in C$ such that
$\rk(A-J)=\rk(A)-2$ and $A-J$ is again semisimple (this is a `semisimple
version' of Proposition \ref{prop:RankRed}). Summarising, this would
prove that for $k \geq 4$ any semisimple element of $R_k$ lies in $kC$.
\end{re}

\subsection{Nilpotent orbits} \label{ssec:Nilp}

As we have seen in Subsection \ref{ssec:2C}, the set
$\overline{2C}\setminus{2C}$ consists of the single nilpotent
orbit $\mathcal{O}[3,2,2,1^{n-7}]$. This motivates the question
of what the minimal $k$ is such that a given nilpotent orbit
$\mathcal{O}$ lies in $kC$. We will see that usually, this $k$ is
just half the rank of $\mathcal{O}$. However, this is not true for
the orbits $\mathcal{O}[3,1^{n-3}]$, $\mathcal{O}[3,2,2,1^{n-3}]$,
and $\mathcal{O}[5,1^{n-5}]$. Furthermore, it remains an open question
what the minimal $k$ is for partitions whose odd entries are all smaller
than $6$. The following lemma will be used to handle odd entries of size
greater than $6$.

\begin{lm}\label{lm:O7}
The nilpotent orbit $\mathcal{O}[7]\subseteq\liea{o}_7$ is contained
in $3C$.
\end{lm}

In the calculation proving this lemma, as in the rest of this subsection,
we compute with concrete matrices that are skew-symmetric with respect
to the bilinear form $(x,y):=\sum_i x_iy_{n+1-i}$.

\begin{proof}
A straightforward computation shows that the difference
\begin{small}
\[
\begin{bmatrix}
0 & 1 & 0 & 0 & 0 & 0 & 0 \\
0 & 0 & 1 & 0 & 0 & 0 & 0 \\ 
0 & 0 & 0 & 1 & 0 & 0 & 0 \\
0 & 0 & 0 & 0 &-1 & 0 & 0 \\ 
0 & 0 & 0 & 0 & 0 &-1 & 0 \\ 
0 & 0 & 0 & 0 & 0 & 0 &-1 \\
0 & 0 & 0 & 0 & 0 & 0 & 0
\end{bmatrix}
-
\begin{bmatrix}
0 & \frac12 & 0 & 1 & 0 & -1 & 0 \\
0 & 0 & 0 & 0 & \frac12 & 0 & 1 \\ 
0 & \frac14 & 0 & \frac12 & 0 & -\frac12 & 0 \\
0 & 0 & 0 & 0 & - \frac12 & 0 &-1 \\ 
0 & 0 & 0 & 0 & 0 & 0 & 0 \\ 
0 & 0 & 0 & 0 & -\frac14 & 0 & -\frac12 \\
0 & 0 & 0 & 0 & 0 & 0 & 0
\end{bmatrix}
\]
\end{small}
\!\!\!\!
is semisimple with three double eigenvalues $0, \frac{i}{2}$,
and $-\frac{i}{2}$, so that it lies in $2C$ by the proof of
Proposition~\ref{prop:2Cbar}. As the matrix on the right lies in
$\mathcal{O}[7]$ and the matrix of the left lies in $C$, this proves
the lemma.
\end{proof}

Recall from the introduction the notation $l(\bd)$ for the number of odd
entries of $\bd$ that are greater than $1$. We are now in a position to
prove Theorem \ref{thm:nilp}.

\begin{proof}[Proof of Theorem \ref{thm:nilp}]
We prove the statement for partitions $[r,r]$ (with $r>1$), $[2r+1,2s+1]$
(with $r>s>0$) and $[2r+1]$ (with $r>2$).  The result will then follow
by pairing equal even entries, pairing odd entries $>1$, taking $[d_1]$
if $l$ is odd, and forming appropriate block matrices.

The orbit $\mathcal{O}[r,r]$ is represented by a $2r \times 2r$-matrix
of the following form (drawn here for $r=3$):
\begin{small}
\[
\left[
\begin{array}{ccc|crr}
0& 1 & &&&\\
&  0 &  1&&&\\
& &  0&&& \\  
\hline
& & &  0 & -1& \\
& & & &  0 & -1 \\
& & & & &  0 
\end{array}  
\right] .    
\]
\end{small}
\!\!\!\!
This matrix is equal to $(E_{1,2}-E_{2r-1,2r})+\dots
+(E_{r-1,r}-E_{r+1,r+2})$ (where $E_{ij}$ is the matrix with an entry 1
at $(i,j)$ and zeroes elsewhere).  These $r-1$ matrices all belong to $C$
(cf. Corollary~\ref{cor:Con}), hence $\mathcal{O}[r,r] \in (r-1)C$.

Next we consider the orbit $\mathcal{O}[2r+1,2s+1]$ with $r>s>0$.
There is a simple recipe for finding a representative of this orbit
\cite[Recipe 5.2.4]{Collingwood93}; by way of example, the orbit
$\mathcal{O}[7,3]$ is represented by (leaving out the off-diagonal
zeroes):
\begin{small}
\[ 
\left[
\begin{array}{ccccc|crrrr}
0&1&&&\\
&0&1&1&\\
&&0&&1&\\
&&&0&1&1\\
&&&&0&&-1\\
\hline
&&&&&0&-1&-1 \\
&&&&&&0&&-1 \\ 
&&&&&&&0&-1\\  
&&&&&&&&0&-1 \\ 
&&&&&&&&&0
\end{array}\right].
\]
\end{small}
\!\!\!
This matrix lies in $4C$, as it is the sum of
$E_{1,2}-E_{9,10},E_{2,3;4}-E_{7;8,9},E_{3;4,5}-E_{6,7;8}$, and
$E_{4,6}-E_{5,7}$, where we use the shorthand notation $E_{i,j_1;j_2}$
for $E_{i,j_1}+E_{i,j_2}$, and its analogue for rows. These matrices all
belong to $C$ by Corollary \ref{cor:Con}, and a moment's reflection shows
that this, too, generalises to the case where $r$ and $s$ are arbitrary,
proving that $\mathcal{O}[2r+1,2s+1] \subseteq (r+s)C$.

Finally, the orbit $\mathcal{O}[2r+1]$ with $r>2$
has a representative of the form drawn in Lemma
\ref{lm:O7} for $r=3$. Subtracting the $(r-3)$ matrices
$E_{1,2}-E_{2r,2r+1},E_{2,3}-E_{2r-1,2r},\ldots,E_{r-3,r-2}-E_{r+4,r+5}
\in C$ yields a matrix with zeroes everywhere except for the $7\times
7$-block of Lemma \ref{lm:O7} in the middle; that lemma shows that this
matrix lies in $3C$. Hence, $\mathcal{O}[2r+1]\subseteq rC$.

To conclude the proof, consider first the case where $l(\bd)$ is even.
We then partition the entries of $\bd$ that are greater than $1$ into
pairs of the forms $[r,r]$ and $[2r+1,2s+1]$ as above. In the case where
$l(\bd)$ is odd and $d_1>5$, we decompose the entries $d_2,\ldots,d_m$
as in the first case, and form the singleton $[d_1]$. In both cases,
a representative of $\mathcal{O}[\bd]$ is then found by gluing 
the block matrices corresponding to the pairs, and in the second case
the block matrix corresponding to the singleton $[d_1]$, together in
an appropriate way. The above calculations show that $\mathcal{O}[\bd]$
lies in $(\rk(A)/2) C$, as claimed.
\end{proof}

We conclude by recalling that, for some nilpotent orbits $\mathcal{O}$,
half the rank of $\mathcal{O}$ does {\em not} suffice!

\begin{lm}
We have $\mathcal{O}[3,1^{n-3}]\subseteq 2C\setminus C$ and
$\mathcal{O}[3,2,2,1^{n-7}], \mathcal{O}[5,1^{n-5}]\subseteq 3C\setminus
2C$.
\end{lm}

\begin{proof}
The first statement follows from Corollary~\ref{cor:Con}: the matrix
corresponding to the partition $[3,1^{n-3}]$ has a non-zero square.
It is clear that the orbits in the second statement are contained in $3C$.
The claim follows then with Lemma~\ref{lm:rk-even}, since both orbits
have rank $4$ and the square of a representative has rank $1$ for the
first, and $3$ for the second orbit.
\end{proof}

The nilpotent orbits of smallest rank for which we do not know
the smallest $k$ such that $kC$ contains them, are therefore
$\mathcal{O}[3,3,3,1^{n-9}]$ and $\mathcal{O}[5,2,2,1^{n-9}]$, which
are both of rank $6$.

\section{Conclusion and further research} \label{sec:Questions}

We have successfully determined, for all classical simple Lie algebras
$\liea{g}$ and all $k \geq 1$, the sets $\overline{kC}$ where $C$
is the adjoint orbit of long root vectors, or, in the terminology of
\cite{Cohen2001}, of extremal elements. If $\liea{g}$ is $\liea{sl}_n$
or $\liea{sp}_{2n}$, then the sets $kC$ are closed, and the minimal $k$
for which they fill the whole space $\liea{g}$ is equal to $n$ or $2n$,
respectively. If, on the other hand, $\liea{g}$ is $\liea{o}_n$, then
$2C$ is not closed, and we only know that the minimal $k$ for which $kC$
is equal to $\liea{o}_n$ lies between $\lfloor \frac{n}{2} \rfloor$
and $\lfloor \frac{n}{2} \rfloor+3$.

We conclude our paper with two rather speculative directions of further
research, suggested by our findings. First, it is shown in \cite{Cohen2001}
that the minimal number of elements of $C$ needed to {\em generate
$\liea{g}$ as a Lie algebra}, is equal to $n$ for $\liea{sl}_n$,
equal to $2n$ for $\liea{sp}_{2n}$, and equal to $\lceil \frac{n}{2}
\rceil$ for $\liea{o}_n$. Of course, the similarity with the numbers
that we listed above may be a coincidence, but if there should be a
direct argument that these numbers are indeed equal, then the results
of \cite{Cohen2001} could be used in solving the remaining open question
concerning $\liea{o}_n$, and in determining the secant varieties of the
minimal orbit in the exceptional Lie algebras, as well.

As mentioned in the introduction, this paper is part of a rather
ambitious project, namely: determining the higher secant varieties of
the minimal orbit in arbitrary irreducible representations of reductive
groups. In that setting, too, the complement of $kC$ in $\overline{kC}$
is worth investigation. The insight that, in the case of $\liea{o}_n$,
the complement of $2C$ in $\overline{2C}$ consists of a nilpotent orbit,
suggests, in the general setting, that $\overline{kC}\setminus kC$ may be
always contained in the null cone. However, if this were true, then it
would follow from our Theorem \ref{thm:nilp} that $\lfloor \frac{n}{2}
\rfloor C$ is already all of $\liea{o}_n$, contrary to what a guess
along the lines of the previous paragraph would yield.


\begin{thebibliography}{10}

\bibitem{Adlandsvik87}
Bj{\o}rn {\AA}dlandsvik.
\newblock {Joins and higher secant varieties.}
\newblock {\em Math. Scand.}, 61(2):213--222, 1987.

\bibitem{Adlandsvik88}
Bj{\o}rn {\AA}dlandsvik.
\newblock {Varieties with an extremal number of degenerate higher secant
  varieties.}
\newblock {\em J. Reine Angew. Math.}, 392:16--26, 1988.

\bibitem{Catalano96}
Michael~L. Catalano-Johnson.
\newblock {The possible dimensions of the higher secant varieties.}
\newblock {\em Am. J. Math.}, 118(2):355--361, 1996.

\bibitem{Catalisano2002}
M.V. Catalisano, A.V. Geramita, and A.~Gimigliano.
\newblock {On the rank of tensors, via secant varieties and fat points.}
\newblock In {\em {Geramita, A. V. (ed.), Zero-dimensional schemes and
  applications. Proceedings of the workshop, Naples, Italy, February 9--12,
  2000. Kingston: Queen's University.}}, volume {123} of {\em {Queen's Pap.
  Pure Appl. Math.}}, pages {135--147}, 2002.

\bibitem{Cohen2003}
Arjeh~M. Cohen and Di\'e~A.H. Gijsbers.
\newblock {\em {\tt GBNP}: a non-commutative {G}r\"{o}bner basis package in
  {\tt GAP}}.
\newblock {\tt http://www.win.tue.nl/\~{ }amc/pub/grobner/doc.html}.

\bibitem{Cohen2001}
Arjeh~M. Cohen, Anja Steinbach, Rosane Ushirobira, and David Wales.
\newblock Lie algebras generated by extremal elements.
\newblock {\em J. Algebra 236}, 236(1):122--154, 2001.

\bibitem{Collingwood93}
David~H. Collingwood and William~M. McGovern.
\newblock {\em Nilpotent orbits in semisimple Lie algebras}.
\newblock Van Nostrand Reinhold Company, New York, 1993.

\bibitem{Fujita81}
Takao Fujita and Joel Roberts.
\newblock {Varieties with small secant varieties: the extremal case.}
\newblock {\em Am. J. Math.}, 103:953--976, 1981.

\bibitem{GAP4}
The GAP~Group.
\newblock {\em {GAP -- Groups, Algorithms, and Programming, Version 4.3}},
  2002.
\newblock {\tt http://www.gap-system.org}.

\bibitem{Goldstein85}
Norman Goldstein.
\newblock {Degenerate secant varieties and a problem on matrices.}
\newblock {\em Pac. J. Math.}, 119:115--124, 1985.

\bibitem{Goodman98}
Roe Goodman and Nolan~R. Wallach.
\newblock {\em Representations and Invariants of the Classical Groups}.
\newblock Cambridge University Press, 1998.

\bibitem{Halmos74}
Paul~R. Halmos.
\newblock {\em {Finite-dimensional vector spaces}}.
\newblock Undergraduate Texts in Mathematics. Springer-Verlag, New York -
  Heidelberg - Berlin, 1974.

\bibitem{Kaji99}
Hajime Kaji, Masahiro Ohno, and Osami Yasukura.
\newblock Adjoint varieties and their secant varieties.
\newblock {\em Indag. Math., New Ser.}, 10(1):45--57, 1999.

\bibitem{Kaji2000}
Hajime Kaji and Osami Yasukura.
\newblock Secant varieties of adjoint varieties: Orbit decomposition.
\newblock {\em J. Algebra}, 227(1):26--44, 2000.

\bibitem{Kraft2001}
Hanspeter Kraft.
\newblock A note on sums of nilpotent matrices of rank one.
\newblock unpublished note, 2001.

\bibitem{Kraft82}
Hanspeter Kraft and Claudio Procesi.
\newblock On the geometry of conjugacy classes in classical groups.
\newblock {\em Comment. Math. Helv.}, 57:539--602, 1982.

\bibitem{Landsberg96}
J.M. Landsberg.
\newblock {On degenerate secant and tangential varieties and local differential
  geometry.}
\newblock {\em Duke Math. J.}, 85(3):605--634, 1996.

\bibitem{Palatini09}
F.~Palatini.
\newblock {Sulle varietà algebriche per le quali sono di dimensione minore
  dell' ordinario, senza riempire lo spazio ambiente, una o alcune delle
  varietà formate da spazi seganti.}
\newblock {\em Torino Atti}, 44:362--375, 1909.

\bibitem{Premet97}
Alexander Premet and Helmut Strade.
\newblock Simple {L}ie algebras of small characteristic. {I}: {S}andwich
  elements.
\newblock {\em J. Algebra}, 189(2):419--480, 1997.

\bibitem{Terracini11}
A.~Terracini.
\newblock {Sulle $V_k$ per cui la varietà degli $S_h(h+1)$-seganti ha
  dimensione minore dell' ordinario.}
\newblock {\em Palermo Rend.}, 31:392--396, 1911.

\bibitem{Zak93}
F.L. Zak.
\newblock {\em {Tangents and secants of algebraic varieties.}}, volume 127 of
  {\em Translations of Mathematical Monographs}.
\newblock American Mathematical Society (AMS), Providence, RI, 1993.

\end{thebibliography}

\end{document}